\documentclass[12pt,twoside]{amsart}
\usepackage{amssymb,amsmath}
\usepackage{a4wide}
\usepackage[utf8]{inputenc}
\usepackage[T1]{fontenc}
\usepackage{times}
\usepackage[hidelinks]{hyperref}
\newcommand{\C}{\mathbb C}
\newcommand{\dbar}{\overline\partial}
\allowdisplaybreaks[2]
\AtBeginDocument{%
  \setlength{\abovedisplayskip}{9pt plus 2pt minus 3pt}%
  \setlength{\belowdisplayskip}{9pt plus 2pt minus 3pt}%
  \setlength{\abovedisplayshortskip}{4pt plus 1pt}%
  \setlength{\belowdisplayshortskip}{6pt plus 1pt minus 2pt}%
}
\numberwithin{equation}{section}
\newtheorem{theorem}{Theorem}[section]
\newtheorem{proposition}[theorem]{Proposition}
\newtheorem{lemma}[theorem]{Lemma}

\theoremstyle{remark}

\begin{document}
\thanks{The author is supported by the Natural Science Foundation of Sichuan Province (No.~2025ZNSFSC0796) and the National Natural Science Foundation of China (No.~12401103).}
\title[The $H^p$-Corona theorem on $C^2$ pseudoconvex domains]
{The $H^p$-Corona theorem on bounded pseudoconvex domains with $C^2$ boundary}
\author[Xu Xing]{Xu Xing}
\begin{abstract}
For $0<p\le2$, we solve the $H^p$-Corona problem on bounded pseudoconvex domains
$\Omega\Subset\C^n$ with $C^2$ boundary, obtaining the bound
$C_\Omega\delta^{-q-1}|\log\delta|$, where
$q=\min\{n,m-1\}$ and $m$ is the number of generators.
The range of $p$ is sharp even on a bounded smooth pseudoconvex domain
admitting a smooth plurisubharmonic defining function and having only one
weakly pseudoconvex boundary point.
\bigskip

\noindent{\textsc{Mathematics Subject Classification} (2020):
32A35, 32W05, 32A10.}

\smallskip

\noindent{\textsc{Keywords}: Corona problem, Hardy spaces,
 $\dbar$-equation.}
\end{abstract}
\address[Xu Xing]{School of Mathematics, Southwest Minzu University,
Chengdu, 610041, China
\newline\mbox{\quad}\,Department of Mathematics and Computer Science,
\newline\mbox{\quad}\,University of Cologne, Cologne, 50931, Germany}
\email{xingxu@swun.edu.cn}
\maketitle

\section{Introduction}

Kakutani posed the Corona problem from a Banach-algebraic viewpoint in 1941
\cite{DouglasHistory}. Equivalently, this problem is whether bounded
holomorphic functions $g_1,\ldots,g_m$ on the unit disc, whose joint
modulus is bounded away from zero, admit bounded holomorphic functions
$h_1,\ldots,h_m$ such that $\sum_{j=1}^m g_jh_j=1$.
Carleson \cite{Carleson} solved the problem in 1962 using what are now known as Carleson measures.
Inspired by H\"ormander's approach, Wolff presented a simplified proof for Carleson's Corona theorem
using boundary $L^\infty$ estimates for the $\dbar$-equation;
see \cite[p.~4]{DouglasHistory}.
The Corona theorem also holds for finitely connected planar domains
and certain infinitely connected domains, including Denjoy domains
\cite{Alling,GarnettJones,Stout}. Nevertheless, the Corona problem remains open for
arbitrary planar domains \cite{DouglasHistory}. By contrast,
Cole constructed an open Riemann surface on which the Corona
theorem fails \cite[Theorem~4.2]{Gamelin}.

Let $\Omega\subset\C^n$ be a domain, and let
$\mathcal F\subset\mathcal O(\Omega)$ be a normed or quasi-normed
function space. Following \cite{ChenXing}, the $\mathcal F$-Corona
problem asks whether, for every
$g=(g_1,\ldots,g_m)\in H^\infty(\Omega)^m$ satisfying
\begin{equation}\label{eq:corona-data}
0<\delta\le |g|:=\Bigl(\sum_{j=1}^m|g_j|^2\Bigr)^{1/2}\le1,
\end{equation}
and every $f\in\mathcal F$, there exists $h:=
(h_1,\ldots,h_m)\in\mathcal F^m$ such that
\[
 g\cdot h:=\sum_{j=1}^m g_jh_j=f,\ \ \ 
 \max_{1\le j\le m}\|h_j\|_{\mathcal F}
 \le C(\delta)\|f\|_{\mathcal F}.
\]
Here $C(\delta)$ may also depend on $\Omega$, $\mathcal F$, and $m$,
but not on $g$ or $f$.
Varopoulos \cite{Varopoulos} obtained a BMOA-Corona theorem on strongly
pseudoconvex domains, subsequently developed by Andersson--Carlsson
\cite{AC}. For $1<p<\infty$, $H^p$-Corona theorems were established
on the unit ball and strongly pseudoconvex domains \cite{Amar91,AC94}.
However, the $H^\infty$-Corona problem remains open even for the
unit ball in $\C^n$, $n\ge2$.
On the other hand, Sibony \cite{Sibony87} and
Forn{\ae}ss--Sibony \cite{FornSibony93} constructed smoothly
bounded pseudoconvex domains in $\C^3$ and $\C^2$, respectively,
on which the $H^\infty$-Corona theorem fails. The domain in
\cite{FornSibony93} is strongly pseudoconvex except at one boundary point.
Andersson \cite{Andersson,Andersson2} proved an $H^2$-Corona
theorem on smooth bounded pseudoconvex domains admitting a
plurisubharmonic defining function, with
$C(\delta)=C_\Omega\delta^{-q-1}|\log\delta|$ for small $\delta>0$,
where $q=\min\{n,m-1\}$.
This result was extended to delta-regular domains in \cite{ChenXing}.
Introduced in \cite{ChenFu}, the class of delta-regular domains includes smooth pseudoconvex
domains of D'Angelo finite type and domains admitting a $C^2$ defining
function that is plurisubharmonic on the boundary.

We establish the following $H^p$-Corona theorem.
\begin{theorem}\label{th:hardy}
Let $\Omega\Subset\C^n$ be a bounded pseudoconvex domain with $C^2$ boundary and let $0<p\le2$. Suppose that $g$ satisfies \eqref{eq:corona-data} for $0<\delta\le e^{-1}$. Then for every $f\in H^p(\Omega)$ there exists $h\in H^p(\Omega)^m$ such that
\begin{equation}\label{eq:main-Hardy}
g\cdot h=f,\ \ \ 
\|h\|_{H^p}\le C_\Omega\delta^{-q-1}|\log\delta|\,\|f\|_{H^p},
\end{equation}
where $q=\min\{n,m-1\}$ and the constant $C_\Omega$ is independent of $p$, $m$, $g$, $f$, and $\delta$.
\end{theorem}

The range of $p$ in Theorem~\ref{th:hardy} is sharp.
In Proposition~\ref{ex:upper-range}, we construct a bounded smooth
pseudoconvex domain in $\C^3$, admitting a smooth plurisubharmonic defining
function, and three bounded holomorphic generators satisfying
\eqref{eq:corona-data} for which $g\cdot h=1$ has no $H^p$ solution
for any $p>2$. Its boundary is strongly pseudoconvex except at one
point of infinite type. Forn{\ae}ss--Sibony \cite{FornSibony93} also showed that $L^p$
estimates for $\dbar$ can fail on smoothly bounded pseudoconvex
domains for $p>2$. 

Following \cite{ChenXing}, we use H\"ormander's $L^2$
method of \cite{Hormander65}, the twisted Morrey--Kohn--H\"ormander
estimate of \cite{Ohsawa,Berndtsson96}, and Skoda's inequality \cite{Skoda}
to prove the main theorem.
The key is Lemma~\ref{lem:normal-row}, which controls
the indefinite Levi term by normal-direction integration by parts,
thereby removing the delta-regularity assumption.

\section{Preliminaries}\label{sec:prelim-Hp}
Throughout the paper, we write
$\partial_\nu=\partial/\partial z_\nu$,
$\partial_{\bar\nu}=\partial/\partial\bar z_\nu$,
$\Delta=4\sum_{\nu=1}^n\partial_\nu\partial_{\bar\nu}$,
and denote by $dV$ the Euclidean volume measure.
Constants denoted by $C_\Omega$ may change from line to line and
depend only on $\Omega$; $a\lesssim b$ means $a\le C_\Omega b$,
and $a\asymp b$ means both $a\lesssim b$ and $b\lesssim a$.

Let $\Omega\Subset\C^n$ be a bounded domain with $C^2$ boundary. For a
$C^2$ defining function $\rho$ and sufficiently small
$t>0$, set
\(
 \Omega_t=\{z\in\Omega:\rho(z)<-t\}.
\)
\subsection{Hardy spaces}
Following \cite[Chapters~I--II]{Stein}, for $0<p<\infty$ define
\begin{equation*}
 H^p(\Omega)=\left\{f\in\mathcal O(\Omega):
 \|f\|_{H^p}^p:=\limsup_{t\downarrow0}
 \int_{\partial\Omega_t}|f|^p\,d\sigma_t
 <\infty\right\},
\end{equation*}
\begin{equation*}
 H^\infty(\Omega)=\left\{f\in\mathcal O(\Omega):
 \|f\|_{H^\infty}:=\sup_{z\in\Omega}|f(z)|<\infty\right\}.
\end{equation*}
For $f=(f_1,\ldots,f_m)\in\mathcal O(\Omega)^m$, write
$|f|^2=\sum_{j=1}^m|f_j|^2$ and set
\[
\|f\|_{H^p}^p:=\limsup_{t\downarrow0}
\int_{\partial\Omega_t}|f|^p\,d\sigma_t.
\]
For $0<p<1$, $\|\cdot\|_{H^p}$ is a quasi-norm.

\begin{lemma}[{cf.~\cite[Chapter~I, \S4]{Stein}, see also \cite{ChenFu}}]
\label{lm:Hardy-limit-Hp}
Let $p>0$ and $h\in\mathcal O(\Omega)^m$. Suppose that $|\nabla\rho|=1$ near
$\partial\Omega$.
Then
\begin{equation}\label{eq:Abel-Hp}
 \lim_{r\uparrow1}(1-r)\int_\Omega|h|^p
 (-\rho)^{-r}\,dV
 =\lim_{t\downarrow0}\int_{\partial\Omega_t}
 |h|^p\,d\sigma_t\ \in[0,\infty].
\end{equation}
If $0<r_j\uparrow1$ and $h_j\in\mathcal O(\Omega)^m$
satisfy
\[
 \liminf_{j\to\infty}(1-r_j)
 \int_\Omega|h_j|^p(-\rho)^{-r_j}\,dV=c_0<\infty,
\]
then a subsequence converges locally uniformly to
$h\in H^p(\Omega)^m$ with $\|h\|_{H^p}^p\le c_0$.
\end{lemma}

\begin{proof}
For $\alpha\ge0$, put
\(
u_\alpha=(|h|^2+\alpha)^{p/2},\
m_\alpha(t)=\int_{\partial\Omega_t}
u_\alpha\,d\sigma_t.
\)
For $\alpha>0$, since $\Delta u_\alpha\ge0$, the
divergence theorem and $|\nabla\rho|=1$ give, for small
$t,s>0$,
\[
\begin{aligned}
m_\alpha(t+s)-m_\alpha(t)
&=\int_{\{t<-\rho<t+s\}}
  \operatorname{div}(u_\alpha\nabla(-\rho))\,dV\\
&=\int_t^{t+s}\int_{\partial\Omega_x}
 \bigl(\nabla u_\alpha\cdot\nabla(-\rho)
       +u_\alpha\Delta(-\rho)\bigr)
 \,d\sigma_x\,dx.
\end{aligned}
\]
Since the outward unit normal to $\Omega_t$ is
$\nabla\rho$, Green's formula yields
\[
\begin{aligned}
m_\alpha'(t)
&=-\int_{\Omega_t}\Delta u_\alpha\,dV
 +\int_{\partial\Omega_t}
   u_\alpha\Delta(-\rho)\,d\sigma_t
\le C_\Omega m_\alpha(t).
\end{aligned}
\]
Consequently,
\[
\frac{d}{dt}\bigl(e^{-C_\Omega t}m_\alpha(t)\bigr)
\le0.
\]
Integrating from $t_1$ to $t_2$, where
$0<t_1<t_2$ are sufficiently small, and letting
$\alpha\downarrow0$, we have
\begin{equation}\label{eq:Hardy-means-Hp}
 m_0(t_2)\le e^{C_\Omega(t_2-t_1)}m_0(t_1),
 \ \ \  m_0(t)=\int_{\partial\Omega_t}|h|^p\,d\sigma_t.
\end{equation}
By the coarea formula, for sufficiently small fixed $t_0>0$,
\[
\begin{aligned}
 (1-r)\int_\Omega|h|^p(-\rho)^{-r}\,dV
 &=(1-r)\int_0^{t_0}x^{-r}m_0(x)\,dx+o(1)
 \ \ \ (r\uparrow1).
\end{aligned}
\]
Since \eqref{eq:Hardy-means-Hp} implies that
$\lim_{t\downarrow0}m_0(t)\in[0,\infty]$ exists and
\[
 (1-r)\int_0^{t_0}x^{-r}\,dx
 =t_0^{1-r}\longrightarrow1,
\]
the Abel limit yields \eqref{eq:Abel-Hp}.
Finally, pass to a subsequence for which
\[
 A_j:=(1-r_j)\int_\Omega
 |h_j|^p(-\rho)^{-r_j}\,dV\longrightarrow c_0.
\]
For $0<x<t$, we have
$m_{0,j}(x)\ge e^{-C_\Omega t}m_{0,j}(t)$ by \eqref{eq:Hardy-means-Hp}, where
$m_{0,j}(t)=\int_{\partial\Omega_t}|h_j|^p\,d\sigma_t$.
Thus
\[
\begin{aligned}
 A_j\ge (1-r_j)\int_0^t x^{-r_j}m_{0,j}(x)\,dx
 \ge e^{-C_\Omega t}t^{1-r_j}m_{0,j}(t).
\end{aligned}
\]
For each fixed sufficiently small $t>0$,
$m_{0,j}(t)$ is uniformly bounded.
Subharmonic estimates and Montel's theorem provide a locally
uniformly convergent subsequence with limit $h\in\mathcal O(\Omega)^m$.
Fatou's lemma yields
$m_0(t)\le e^{C_\Omega t}c_0$; letting $t\downarrow0$
proves $\|h\|_{H^p}^p\le c_0$.
\end{proof}

\subsection{Carleson-type inequalities}

Following \cite[Lemmas~2.1--2.2]{ChenXing}, we record the
following estimates.
\begin{lemma}\label{lm:data-Green-Hp}
Let $u\in C^2(\Omega)$ be positive with $\log u$
subharmonic, and let $\psi\in C^2(\Omega)$ be bounded and
subharmonic. For sufficiently small $t>0$,
\begin{align}
 \int_{\Omega_t}(-\rho-t)\Delta u\,dV
 &\le C_\Omega\int_{\partial\Omega_t}
 u\,d\sigma_t,
 \label{eq:positive-Green}\\
 \int_{\Omega_t}(-\rho-t)u
 \bigl(\|\psi\|_\infty\Delta\psi+|\partial\psi|^2\bigr)\,dV
 &\le C_\Omega\|\psi\|_\infty^2
 \int_{\partial\Omega_t}u\,d\sigma_t.
 \label{eq:weighted-Green-Hp}
\end{align}
The constant is independent of $u$, $\psi$, and $t$.
\end{lemma}

\begin{proof}
By the construction in \cite[Lemma~2.1]{ChenXing}, choose a
negative $C^2$ defining function $\varrho$ with
$\Delta\varrho\ge c_\Omega>0$ on $\Omega$ and
$\varrho=\rho+B\rho^2$ near $\partial\Omega$,
for some fixed $B>0$.
Set $\varrho_t=\varrho+t-Bt^2$.
For every nonnegative subharmonic $\varphi\in C^2(\Omega)$,
Green's formula yields
\[
\begin{aligned}
\int_{\Omega_t}
(-\varrho_t)\Delta \varphi\,dV
+\int_{\Omega_t}
\varphi\Delta\varrho\,dV
=
\int_{\partial\Omega_t}
\varphi\,\partial_{\mathbf n}\varrho_t
\,d\sigma_t.
\end{aligned}
\]
Since $-\varrho_t\asymp-\rho-t$ on $\Omega_t$
and $\Delta\varrho\ge c_\Omega>0$, we have
\[
 \int_{\Omega_t}\varphi\,dV
 +\int_{\Omega_t}(-\rho-t)\Delta \varphi\,dV
 \le C_\Omega\int_{\partial\Omega_t}
 \varphi\,d\sigma_t.
\]
Taking $\varphi=u$ proves \eqref{eq:positive-Green}.
For $M=\|\psi\|_\infty>0$, put $s=(\psi+M)/(2M)$.
Since $\log u$ is subharmonic and $0\le s\le1$, applying the
same inequality to $ue^s$ and $ue^{s^2}$ yields
\eqref{eq:weighted-Green-Hp}. The case $\|\psi\|_\infty=0$ is immediate.
\end{proof}

\subsection{Basic $\dbar$ estimates}
\label{subsec:constrained-twist}

Following \cite{HormanderBook,Ohsawa,ChenXing}, let
$v=\sum_{\nu=1}^n v_\nu\,d\bar z_\nu$ and
$\omega=\sum_{\nu=1}^n\omega_\nu\,d\bar z_\nu$
be $\mathbb C^m$-valued $(0,1)$-forms, with
$v_\nu=(v_{\nu,1},\ldots,v_{\nu,m})$ and
$\omega_\nu=(\omega_{\nu,1},\ldots,\omega_{\nu,m})$.
For a $\mathbb C^m$-valued function $g=(g_1,\ldots,g_m)$
and a function $\psi$, we use
\[
\begin{aligned}
\langle v,\omega\rangle
=\sum_{\nu=1}^n\sum_{\mu=1}^m
v_{\nu,\mu}\overline{\omega_{\nu,\mu}},
\ \ \ 
\partial\psi\mathbin{\lrcorner}\omega
=\sum_{\nu=1}^n(\partial_\nu\psi)\omega_\nu
\in\mathbb C^m,\ \ \ 
\partial g\mathbin{\lrcorner}\omega
=\sum_{\nu=1}^n\sum_{\mu=1}^m
(\partial_\nu g_\mu)\omega_{\nu,\mu}\in\mathbb C.
\end{aligned}
\]
For a real $(1,1)$-form
$\Theta=i\sum_{\mu,\nu}\Theta_{\mu\bar\nu}
 dz_\mu\wedge d\bar z_\nu$, we denote
\begin{equation*}
 \langle[\Theta,\Lambda]\omega,\omega\rangle
 :=\sum_{\mu,\nu=1}^n
 \Theta_{\mu\bar\nu}\langle\omega_\mu,\omega_\nu\rangle,
 \ \ \ 
 \langle[\Theta,\Lambda]^{-1}v,v\rangle
 :=\sup_{\langle[\Theta,\Lambda]\omega,\omega\rangle\le1}
 |\langle v,\omega\rangle|^2.
\end{equation*}
The inverse form is understood when $\Theta\ge0$.
We also write
\[
 |\bar\nabla\omega|^2
 =\sum_{\mu,\nu}|\partial_{\bar\nu}\omega_\mu|^2.
\]
Let $D\Subset\Omega$ be a pseudoconvex domain with $C^2$ boundary and a $C^2$ defining
function $\rho_D$.
For $\varphi\in C^2(\overline D)$ and $\eta>0$, put
\[
 (v,\omega)_\varphi=\int_D\langle v,\omega\rangle
 e^{-\varphi}\,dV,
 \ \ \ 
 \|\omega\|_{\eta,\varphi}^2
 =\int_D\eta|\omega|^2e^{-\varphi}\,dV.
\]
Let $T$ and $S$ be the closed, densely defined operators
induced by $\dbar$,
\[
 T:L^2(D,e^{-\varphi})^m
   \longrightarrow L^2_{(0,1)}(D,e^{-\varphi})^m,
 \ \ \ 
 S:L^2_{(0,1)}(D,e^{-\varphi})^m
   \longrightarrow L^2_{(0,2)}(D,e^{-\varphi})^m.
\]
Let $T_\varphi^*$ denote the weighted adjoint of $T$. For a smooth $(0,1)$-form $\omega\in\operatorname{Dom}T_\varphi^*$, 
\[
 T^*_\varphi\omega
 =-\sum_{\nu=1}^n(\partial_\nu-\varphi_\nu)\omega_\nu.
\]

Assume that
$g\in\mathcal O(\Omega)^m$ satisfies \eqref{eq:corona-data}.
Define
\[
 \mathcal H_j=\{\omega\in L^2_{(0,j)}(D,e^{-\varphi})^m:
                     g\cdot\omega=0\},\quad j=0,1,2,
\]
with $T_{\mathcal H}=\dbar:\mathcal H_0\to\mathcal H_1$,
$S_{\mathcal H}=\dbar:\mathcal H_1\to\mathcal H_2$. Set
\[
 \mathcal D_1(D)=\left\{\omega\in C^1_{(0,1)}(\overline D)^m:
 g\cdot\omega=0,\
 \partial\rho_D\mathbin{\lrcorner}\omega=0
 \ \text{on }\partial D\right\}.
\]
Write
\(
 \mathcal G_\varphi
 :=\operatorname{Dom}T^*_{\mathcal H,\varphi}
   \cap\operatorname{Dom}S_{\mathcal H}.
\)
For $\omega\in\mathcal H_1$, set
\begin{equation}\label{eq:Phi-g}
 \Phi_g\omega
 :=\frac{\bar g}{|g|^2}
 (\partial g\mathbin{\lrcorner}\omega).
\end{equation}
Following \cite{Ohsawa} and
\cite[Lemmas~4.1--4.2]{ChenXing}, we have,
for $\omega\in\mathcal D_1(D)$,
\[
 T^*_{\mathcal H,\varphi}\omega
 =T^*_{\varphi}\omega-\Phi_g\omega,
 \ \ \ 
 \langle T^*_{\mathcal H,\varphi}\omega,\Phi_g\omega\rangle=0.
\]
Moreover, $\mathcal D_1(D)$ is dense in
$\mathcal G_\varphi$ in the graph norm.

\begin{lemma}[
{\cite[(4.4)]{ChenXing}}]\label{lem:KMH-basic}
For $0<r<1$, set $\eta=(-\rho)^r$. Then for
 every $\omega\in\mathcal D_1(D)$,
\begin{align}\label{eq:KMH-basic}
& \frac1{1-r}\|T^*_{\mathcal H,\varphi}\omega\|_{\eta,\varphi}^2
  +\|S_{\mathcal H}\omega\|_{\eta,\varphi}^2
  +\|\Phi_g\omega\|_{\eta,\varphi}^2\notag\\
&\ge\int_D\eta\bigl(
       \langle[i\partial\bar\partial\varphi,\Lambda]
                         \omega,\omega\rangle
          +|\bar\nabla\omega|^2\bigr)e^{-\varphi}\,dV
+r\int_D(-\rho)^{r-1}
      \langle[i\partial\bar\partial\rho,\Lambda]
                       \omega,\omega\rangle e^{-\varphi}\,dV.
\end{align}
\end{lemma}
The proof of this twisted Morrey--Kohn--H\"ormander inequality
is inspired by the twisted estimates of Ohsawa \cite{Ohsawa}
and Berndtsson \cite{Berndtsson96}.

\section{Proof of Theorem~\ref{th:hardy}}
\label{sec:Hardy-proof}

From now on, assume that $\Omega$ is pseudoconvex. Let
$d(z)=\operatorname{dist}(z,\partial\Omega)$ in $\Omega$.
Choose $A>0$, to be fixed below, and
$0<\varepsilon_1<\varepsilon_2$ sufficiently small. Put
$U=\{z\in\Omega:0<d(z)<\varepsilon_2\}$ and
$V=\{z\in\Omega:0<d(z)<\varepsilon_1\}$, so that
\[
 d\in C^2(U),\ \ \ |\nabla d|=1,\ \ \ 
 |\nabla^2d|\le C_\Omega,\ \ \ 1+2Ad>0\quad\text{on }U.
\]
Choose a defining function $\rho\in C^2(\overline\Omega)$,
with $-\rho\asymp d$ and
\begin{equation}\label{eq:rho-choice}
 \rho=-d-Ad^2\quad\text{on }U.
\end{equation}
Fix $0<\varepsilon_0<\varepsilon_1$ such that
$\varepsilon_0+A\varepsilon_0^2<\min_{\{d\ge\varepsilon_2\}}(-\rho)$.
For $0<\varepsilon<\varepsilon_0$, the function
$\rho_\varepsilon:=\rho+\varepsilon+A\varepsilon^2$ defines
$\Omega_\varepsilon$, since on $U$,
$\rho_\varepsilon=-(d-\varepsilon)(1+A(d+\varepsilon))$.
By Oka's lemma, $\Omega_\varepsilon$ is
pseudoconvex with $C^2$ boundary.

Choose $\vartheta\in C^\infty(\overline\Omega)$ with
$0\le\vartheta\le1$, $\vartheta=1$ on $V$, and
$\vartheta=0$ on $\{d\ge\varepsilon_2\}$.
For a $\C^m$-valued $(0,1)$-form
$\omega$, define on $U$
\begin{equation}\label{eq:normal-row-fields}
 \begin{gathered}
 L=\sum_{\lambda=1}^n L_\lambda\partial_\lambda:=\frac{\vartheta}{|\partial\rho|^2}
   \sum_{\lambda=1}^n\rho_{\bar\lambda}\partial_\lambda,
 \ \ \ 
  b=-\bar L(\partial\rho)\mathbin{\lrcorner}\omega=-\frac{\vartheta}{|\partial\rho|^2}
   \left(\sum_{\lambda=1}^n
   \rho_\lambda\partial\rho_{\bar\lambda}\right)
   \mathbin{\lrcorner}\omega.
 \end{gathered}
\end{equation}
Set $L=b=0$ on $\Omega\setminus U$. Then $L\in C^1$,
$|b|\le C_\Omega|\omega|$, and $L\rho=\vartheta$ on $\Omega$.

\begin{lemma}\label{lem:normal-row}
For $A=A(\Omega)>0$ sufficiently large, every $\C^m$-valued $(0,1)$-form $\omega$ satisfies
on $\Omega$
\begin{equation}\label{eq:normal-row-lower}
 \langle[i\partial\bar\partial\rho,\Lambda]\omega,\omega\rangle
 \ge 2\operatorname{Re}\langle
   \partial(-\rho)\mathbin{\lrcorner}\omega,b\rangle
      -C_\Omega(-\rho)|\omega|^2.
\end{equation}
For every $0<\varepsilon<\varepsilon_0$,
real-valued $\eta,\varphi\in C^2(\overline{\Omega_\varepsilon})$
with $\eta>0$,
and $\omega\in\mathcal D_1(\Omega_\varepsilon)$,
\begin{equation}\label{eq:normal-row-adjoint-bound}
\begin{aligned}
2\operatorname{Re}\int_{\Omega_\varepsilon}
  \langle\partial\eta\mathbin{\lrcorner}\omega,b\rangle
                                  e^{-\varphi}\,dV
 \ge-\|T^*_{\mathcal H,\varphi}\omega\|_{\eta,\varphi}^2
       -\tfrac12\|\bar\nabla\omega\|_{\eta,\varphi}^2
       -C_0\|\omega\|_{\eta,\varphi}^2.
\end{aligned}
\end{equation}
The constant $C_0=C_0(\Omega)$ is independent of $\eta,\varphi,\varepsilon,g,m$.
\end{lemma}

\begin{proof}
Oka's lemma implies, on $U\cap\ker\partial\rho$,
\[
i\partial\bar\partial\rho
=d(1+2Ad)\,i\partial\bar\partial(-\log d)\ge0.
\]
 Moreover, where $\vartheta=1$,
\[
 \sum_{\mu,\nu}\rho_{\mu\bar\nu}L_\mu\overline{L_\nu}
 =-\frac{1}{1+2Ad}\sum_{\mu,\nu}
 d_{\mu\bar\nu}\frac{d_{\bar\mu}d_\nu}{|\partial d|^4}
 -\frac{2A}{(1+2Ad)^2}\le0
\]
for $A$ sufficiently large and $U$ sufficiently small.

On $V$, $L\rho=1$. Thus
$\omega^\top_\mu=\omega_\mu+
L_\mu(\partial(-\rho)\mathbin{\lrcorner}\omega)$ satisfies
$\partial\rho\mathbin{\lrcorner}\omega^\top=0$, and
\[
\begin{aligned}
\langle[i\partial\bar\partial\rho,\Lambda]\omega,\omega\rangle
&=\langle[i\partial\bar\partial\rho,\Lambda]
\omega^\top,\omega^\top\rangle
+2\operatorname{Re}\langle
\partial(-\rho)\mathbin{\lrcorner}\omega,b\rangle
-\left(\sum_{\mu,\nu}\rho_{\mu\bar\nu}
L_\mu\overline{L_\nu}\right)
|\partial\rho\mathbin{\lrcorner}\omega|^2\\
&\ge 2\operatorname{Re}\langle
\partial(-\rho)\mathbin{\lrcorner}\omega,b\rangle.
\end{aligned}
\]

On $\Omega\setminus V$, we have $-\rho\ge c_\Omega>0$.
Since $|\partial\rho|+|\partial\bar\partial\rho|\le C_\Omega$
and $|b|\le C_\Omega|\omega|$,
\[
\begin{aligned}
\langle[i\partial\bar\partial\rho,\Lambda]\omega,\omega\rangle
-2\operatorname{Re}
\langle\partial(-\rho)\mathbin{\lrcorner}\omega,b\rangle
\ge -C_\Omega|\omega|^2
\ge -C_\Omega(-\rho)|\omega|^2.
\end{aligned}
\]
This proves \eqref{eq:normal-row-lower}.

Write $x_1,\ldots,x_{2n}$ for real coordinates.
Differentiating $|\nabla d|^2=1$ first in $z_\mu$
and then in $\bar z_\nu$ yields, in distributions,
\[
 \sum_{j=1}^{2n}d_{x_j}d_{\mu x_j}=0,\ \ \ 
 (\nabla d\cdot\nabla)d_{\mu\bar\nu}
 =-\sum_{j=1}^{2n}d_{\mu x_j}d_{\bar\nu x_j}.
\]
Since $\rho_{\mu\bar\nu}=-(1+2Ad)d_{\mu\bar\nu}-2Ad_\mu d_{\bar\nu}$ and
$L+\bar L=-2\vartheta(1+2Ad)^{-1}\nabla d\cdot\nabla$,
\[
 (L+\bar L)\rho_{\mu\bar\nu}
 =\frac{4A\vartheta}{1+2Ad}d_{\mu\bar\nu}
   -2\vartheta\sum_{j=1}^{2n}d_{\mu x_j}d_{\bar\nu x_j},
 \ \ \ 
 \|(L+\bar L)\rho_{\mu\bar\nu}\|_\infty\le C_\Omega.
\]
The distributional identity
\(
 \partial_\mu\rho_{\lambda\bar\nu}
 =\partial_\lambda\rho_{\mu\bar\nu}
\)
shows that, with summation over $\lambda$,
\[
 \partial_\mu(L_\lambda\rho_{\lambda\bar\nu})
 =(\partial_\mu L_\lambda)\rho_{\lambda\bar\nu}
 +L_\lambda\partial_\lambda\rho_{\mu\bar\nu}.
\]
For  $\omega\in\mathcal D_1(\Omega_\varepsilon)$, taking the Hermitian part yields
\begin{equation}\label{eq:real-normal-hessian}
\begin{aligned}
\operatorname{Re}\sum_{\mu,\nu,\lambda}
 \partial_\mu(L_\lambda\rho_{\lambda\bar\nu})
 \langle\omega_\mu,\omega_\nu\rangle
=\frac12\sum_{\mu,\nu}
 ((L+\bar L)\rho_{\mu\bar\nu})
 \langle\omega_\mu,\omega_\nu\rangle
 +O_\Omega(|\omega|^2).
\end{aligned}
\end{equation}
Thus the Hermitian part of these coefficient derivatives is
represented by a bounded continuous matrix.
Fix $0<\varepsilon<\varepsilon_0$ and take a standard real
mollifier $\zeta_s$ supported in $\{|x|<s\}$, where
$0<s<\varepsilon/2$. Since
$\operatorname{dist}(\overline{\Omega_\varepsilon},
\partial\Omega)=\varepsilon$, all convolutions below are
well-defined on $\overline{\Omega_\varepsilon}$. Put
\[
 b_s=-\left[
 \left(\frac{\vartheta}{|\partial\rho|^2}
 \sum_{\lambda=1}^n\rho_\lambda\partial\rho_{\bar\lambda}
 \right)*\zeta_s\right]\mathbin{\lrcorner}\omega.
\]
Then $g\cdot b_s=0$.
By
\eqref{eq:real-normal-hessian}, $b_s$ converge uniformly on
$\overline{\Omega_\varepsilon}$ and are bounded independently
of $s$ and $\varepsilon$.
Since
$\partial\rho_\varepsilon\mathbin{\lrcorner}\omega=0$
on $\partial\Omega_\varepsilon$,
letting $s\downarrow0$ and using
\eqref{eq:real-normal-hessian} gives
\[
\begin{aligned}
\operatorname{Re}\int_{\Omega_\varepsilon}
 \langle\partial\eta\mathbin{\lrcorner}\omega,b\rangle
 e^{-\varphi}\,dV
&=\operatorname{Re}\int_{\Omega_\varepsilon}\eta
 \Bigl(\langle T_{\mathcal H,\varphi}^*\omega,b\rangle
 +\sum_{\mu,\nu,\lambda}L_\lambda\rho_{\lambda\bar\nu}
 \langle\omega_\mu,\partial_{\bar\mu}\omega_\nu\rangle
 \Bigr)e^{-\varphi}\,dV\\
&\quad+\frac12\int_{\Omega_\varepsilon}\eta
 \sum_{\mu,\nu}((L+\bar L)\rho_{\mu\bar\nu})
 \langle\omega_\mu,\omega_\nu\rangle e^{-\varphi}\,dV
 +O_\Omega(\|\omega\|_{\eta,\varphi}^2).
\end{aligned}
\]
Since $|b|\le C_\Omega|\omega|$ and
$\|(L+\bar L)\rho_{\mu\bar\nu}\|_\infty\le C_\Omega$,
Cauchy--Schwarz gives
\[
\begin{aligned}
 \left|2\operatorname{Re}\int_{\Omega_\varepsilon}
     \langle\partial\eta\mathbin{\lrcorner}\omega,b\rangle
                                     e^{-\varphi}\,dV\right|
 \le 2C_1\|\omega\|_{\eta,\varphi}
 \bigl(\|T^*_{\mathcal H,\varphi}\omega\|_{\eta,\varphi}
       +\|\bar\nabla\omega\|_{\eta,\varphi}
       +\|\omega\|_{\eta,\varphi}\bigr).
\end{aligned}
\]
Young's inequality $2ab\le\tau a^2+\tau^{-1}b^2$ yields
\[
\begin{aligned}
2C_1\|\omega\|_{\eta,\varphi}
 \|T^*_{\mathcal H,\varphi}\omega\|_{\eta,\varphi}
&\le\|T^*_{\mathcal H,\varphi}\omega\|_{\eta,\varphi}^2
 +C_1^2\|\omega\|_{\eta,\varphi}^2,\\
2C_1\|\omega\|_{\eta,\varphi}
 \|\bar\nabla\omega\|_{\eta,\varphi}
&\le\tfrac12\|\bar\nabla\omega\|_{\eta,\varphi}^2
 +2C_1^2\|\omega\|_{\eta,\varphi}^2.
\end{aligned}
\]
Thus \eqref{eq:normal-row-adjoint-bound} holds with
$C_0=3C_1^2+2C_1$.
\end{proof}

\begin{proposition}\label{prop:profile-basic}
There exists $C_\Omega>0$ such that, for
$0<\varepsilon<\varepsilon_0$, any real
$\varphi\in C^2(\overline{\Omega_\varepsilon})$,
$0\le r<1$, and $\omega\in\mathcal D_1(\Omega_\varepsilon)$,
with $\eta=(-\rho)^r$,
\begin{align}\label{eq:profile-basic}
 &\int_{\Omega_\varepsilon}\eta\left(
     \langle[i\partial\bar\partial\varphi,\Lambda]
                         \omega,\omega\rangle
     -|\Phi_g\omega|^2-C_\Omega|\omega|^2
     +\tfrac12|\bar\nabla\omega|^2\right)e^{-\varphi}\,dV\notag\\
 &\ \ \ \le\frac{2-r}{1-r}
     \|T^*_{\mathcal H,\varphi}\omega\|_{\eta,\varphi}^2
     +\|S_{\mathcal H}\omega\|_{\eta,\varphi}^2.
\end{align}
The estimate extends to $\mathcal G_\varphi$.
\end{proposition}

\begin{proof}
First assume $0<r<1$.
Since
\(
 \partial\eta=r(-\rho)^{r-1}\partial(-\rho),
\)
\eqref{eq:normal-row-lower} and
\eqref{eq:normal-row-adjoint-bound} imply
\[
\begin{aligned}
r\int_{\Omega_\varepsilon}(-\rho)^{r-1}
 \langle[i\partial\bar\partial\rho,\Lambda]
                      \omega,\omega\rangle e^{-\varphi}\,dV
&\ge
 2\operatorname{Re}\int_{\Omega_\varepsilon}
   \langle\partial\eta\mathbin{\lrcorner}\omega,b\rangle
                   e^{-\varphi}\,dV
   -C_\Omega\|\omega\|_{\eta,\varphi}^2\\
&\ge-\|T^*_{\mathcal H,\varphi}\omega\|_{\eta,\varphi}^2
 -\tfrac12\|\bar\nabla\omega\|_{\eta,\varphi}^2
 -C_\Omega\|\omega\|_{\eta,\varphi}^2.
\end{aligned}
\]
Then \eqref{eq:KMH-basic}
 yields
\eqref{eq:profile-basic}. For $r=0$, the claim follows
from \eqref{eq:KMH-basic} by taking $r\downarrow0$.
Graph-norm approximation and weak lower semicontinuity of
$\bar\nabla\omega$ yield the extension to $\mathcal G_\varphi$.
\end{proof}

\begin{lemma}[Skoda's inequality {\cite[Lemma~4.3]{ChenXing}}]\label{lem:Skoda}
For any $\C^m$-valued $(0,1)$-tensor $\omega$,
\begin{equation}\label{eq:Skoda-tensor}
 |g|^2\left|\left\langle
 \dbar\left(\frac{\bar g}{|g|^2}\right),\omega
 \right\rangle\right|^2
 \le q\langle[\Theta_g,\Lambda]\omega,\omega\rangle,
\end{equation}
where $\Theta_g=i\partial\bar\partial\log|g|^2$. Consequently, for $f\in\mathcal O(\Omega)$,
\begin{equation}\label{eq:Skoda-source-general}
 \left\langle[\Theta_g,\Lambda]^{-1}
 \left(f\dbar\frac{\bar g}{|g|^2}\right),
 f\dbar\frac{\bar g}{|g|^2}\right\rangle
 \le q\frac{|f|^2}{|g|^2}.
\end{equation}
\end{lemma}
\begin{proof}[Proof of Theorem~\ref{th:hardy}]
The cases $f\equiv0$ and $m=1$ are immediate. Put
\(
 v=f\dbar\,({\bar g}/{|g|^2}).
\)
Then $g\cdot v=0$ and $S_{\mathcal H}v=0$. It suffices to find
$u\in\mathcal H_0$ with $\dbar u=v$ and set
\[
 h=\frac{f\bar g}{|g|^2}-u,
 \ \ \  \dbar h=0,\quad g\cdot h=f.
\]

Fix $0<p\le2$ and $\alpha>0$, and put
\(s=p/2,\ \gamma=|\log\delta|^{-1},\
 \chi_\alpha=(1-s)\log(|f|^2+\alpha)\).
\begin{equation}\label{eq:common-weights}
\begin{gathered}
 \varphi_\alpha=M|z|^2+(q+\gamma)\log|g|^2+\chi_\alpha,
 \ \ \ \Theta=i\partial\bar\partial|z|^2+\gamma\Theta_g.
\end{gathered}
\end{equation}
Here $M= C_\Omega+1$, where $C_\Omega$ is the constant
in Proposition~\ref{prop:profile-basic}. 
Since $g\cdot\omega=0$, we have
\[
\left\langle
\dbar\left(\frac{\bar g}{|g|^2}\right),\omega
\right\rangle
=\frac{\overline{\partial g\mathbin{\lrcorner}\omega}}{|g|^2}.
\]
Hence, by Lemma~\ref{lem:Skoda},
\[
\begin{aligned}
|\Phi_g\omega|^2
=\frac{|\partial g\mathbin{\lrcorner}\omega|^2}{|g|^2}
=|g|^2\left|
\left\langle
\dbar\left(\frac{\bar g}{|g|^2}\right),\omega
\right\rangle\right|^2
\le q\langle[\Theta_g,\Lambda]\omega,\omega\rangle.
\end{aligned}
\]
For $0<\varepsilon<\varepsilon_0$ and
$\omega\in\mathcal D_1(\Omega_\varepsilon)$,
\[
\begin{aligned}
 \langle[i\partial\bar\partial\varphi_\alpha,\Lambda]
           \omega,\omega\rangle-|\Phi_g\omega|^2
                                      -C_\Omega|\omega|^2
\ge \langle[\Theta+i\partial\bar\partial\chi_\alpha,
                              \Lambda]\omega,\omega\rangle.
\end{aligned}
\]
Therefore, Proposition~\ref{prop:profile-basic} yields
\begin{align}\label{eq:profile-applied}
 &\int_{\Omega_\varepsilon}(-\rho)^r
  \Bigl(\langle[\Theta+i\partial\bar\partial\chi_\alpha,
                 \Lambda]\omega,\omega\rangle
             +\tfrac12|\bar\nabla\omega|^2\Bigr)
                                      e^{-\varphi_\alpha}\,dV\notag\\
 &\quad\le\frac{2-r}{1-r}
     \|T^*_{\mathcal H,\varphi_\alpha}\omega\|_{(-\rho)^r,\varphi_\alpha}^2
      +\|S_{\mathcal H}\omega\|_{(-\rho)^r,\varphi_\alpha}^2,
 \ \ \  0\le r<1.
\end{align}

For $1\le\lambda\le n$, we have
\(
 \partial_\lambda\log|g|^2
 ={\langle\partial_\lambda g,g\rangle}/{|g|^2}.
\)
Consequently,
\[
\begin{aligned}
\sum_{\lambda=1}^n L_\lambda
\bar\partial\partial_\lambda\log|g|^2
=\sum_{\lambda=1}^n\sum_{k=1}^m
L_\lambda(\partial_\lambda g_k)
\bar\partial\left(\frac{\bar g_k}{|g|^2}\right)
=(Lg)\cdot\bar\partial\left(\frac{\bar g}{|g|^2}\right).
\end{aligned}
\]
Thus, componentwise,
\begin{align}\label{eq:vg-source}
v_{g,\alpha}
&:=Lv+\frac{\bar g}{|g|^2}(Lg)\cdot v
  -(L\varphi_\alpha)v\notag\\
&=\left(Lf-fL\log|g|^2
  -fL\varphi_\alpha\right)
  \dbar\left(\frac{\bar g}{|g|^2}\right)\notag\\
&=\left\{Lf-fL\chi_\alpha
  -fL\left(M|z|^2+
  (q+1+\gamma)\log|g|^2\right)\right\}
  \dbar\left(\frac{\bar g}{|g|^2}\right).
\end{align}
The second term in the first line is orthogonal to
$\omega\in\mathcal H_1$.
Differentiating $g\cdot v=0$ along $L$ also shows that
$g\cdot v_{g,\alpha}=0$.
Combining $\Theta\ge\gamma\Theta_g$ with Lemma~\ref{lem:Skoda} and the Cauchy--Schwarz inequality yields
\[
\begin{gathered}
 \left\langle[\Theta,\Lambda]^{-1}
    \dbar\!\left(\frac{\bar g}{|g|^2}\right),
    \dbar\!\left(\frac{\bar g}{|g|^2}\right)\right\rangle
       \le\frac{q}{\gamma|g|^2},
\end{gathered}
\]
where \( |\partial g|^2=\sum_\nu|\partial_\nu g|^2\). Moreover, as $\delta\le|g|\le1$ and
$\chi_\alpha=(1-s)\log(|f|^2+\alpha)$,
\[
 |g|^{-2\gamma}\le e^2,\ \ \ 
 e^{-\varphi_\alpha+\chi_\alpha}\le e^2\delta^{-2q},
 \ \ \  |f|^2e^{-\chi_\alpha}\le|f|^p.
\]
Since $L$ is bounded, $q\le n$, and $|\log\delta|\ge1$,
we obtain the following estimate almost everywhere:

\begin{align}\label{eq:unified-source-bound}
&\left(|v|^2+\langle[\Theta,\Lambda]^{-1}
v_{g,\alpha},v_{g,\alpha}\rangle\right)e^{-\varphi_\alpha}
\notag\\
&\le C_\Omega e^{-\varphi_\alpha}
\left[
\frac{|f|^2|\partial g|^2}{|g|^4}
+\frac{q}{\gamma|g|^2}
\left|Lf-fL\chi_\alpha
-fL\left(M|z|^2+
(q+1+\gamma)\log|g|^2\right)\right|^2
\right]\notag\\
&\le \frac{C_\Omega q}{\gamma}e^{-\varphi_\alpha}
\left[
\frac{|\partial f-f\partial\chi_\alpha|^2
+|f|^2}{|g|^2}
+(q+2)^2|f|^2\frac{|\partial g|^2}{|g|^4}
\right]\notag\\
&\le C_\Omega\delta^{-2q}|\log\delta|
\left[
\delta^{-2}|\partial f-f\partial\chi_\alpha|^2
e^{-\chi_\alpha}
+\delta^{-2}|f|^p
+|f|^p\frac{|\partial g|^2}{|g|^4}
\right].
\end{align}
Put $u_\alpha=(|f|^2+\alpha)^s$. Direct differentiation yields
\[
\begin{aligned}
 \partial f-f\partial\chi_\alpha
 =\frac{\alpha+s|f|^2}{|f|^2+\alpha}\partial f,\ \ \ 
 \Delta u_\alpha
 =4s(|f|^2+\alpha)^{s-2}(\alpha+s|f|^2)|\partial f|^2.
\end{aligned}
\]
On $\{f\ne0\}$,
\[
 4s^2(|f|^2+\alpha)^{s-1}|\partial f|^2
 \le\Delta u_\alpha
 \le 4s|f|^{2s-2}|\partial f|^2.
\]
For $0\le\zeta\in C_c^\infty(\Omega)$,
\[
\begin{aligned}
4s^2\int_\Omega\zeta(|f|^2+\alpha)^{s-1}
|\partial f|^2\,dV
\le\int_\Omega u_\alpha\Delta\zeta\,dV
\longrightarrow\int_\Omega |f|^{2s}\Delta\zeta\,dV.
\end{aligned}
\]
Monotone convergence in the preceding inequality shows that
$|f|^{p-2}|\partial f|^2\in L^1_{\mathrm{loc}}(\Omega)$.
The upper bound for $\Delta u_\alpha$ and dominated convergence then imply
\[
 \Delta|f|^p=p^2|f|^{p-2}|\partial f|^2
 \quad\text{in distributions}.
\]
For $\psi=\log((|g|^2+\delta^2)/2)$,
\[
 \|\psi\|_\infty\le2|\log\delta|,\ \ \ 
 \Delta\psi\ge\frac{4\delta^2|\partial g|^2}
 {(|g|^2+\delta^2)^2}
 \ge\delta^2\frac{|\partial g|^2}{|g|^4}.
\]
Applying Lemma~\ref{lm:data-Green-Hp} to $u_\alpha$ and $\psi$,
using also
\[
 \int_{\Omega_\varepsilon}u_\alpha\,dV
 \le C_\Omega\int_{\partial\Omega_\varepsilon}
 u_\alpha\,d\sigma_\varepsilon
\]
from its proof, and letting first $\alpha\downarrow0$ and then
$\varepsilon\downarrow0$, we obtain
\[
\begin{aligned}
 \int_\Omega d\,\Delta|f|^p\,dV
 &\le C_\Omega\|f\|_{H^p}^p,\\
 \int_\Omega d\,|f|^p\frac{|\partial g|^2}{|g|^4}\,dV
 &\le C_\Omega\delta^{-2}|\log\delta|\,\|f\|_{H^p}^p,\\
 \int_\Omega d\,|f|^p\,dV
 &\le C_\Omega\|f\|_{H^p}^p.
\end{aligned}
\]
Since, almost everywhere,
\[
\begin{aligned}
 |\partial f-f\partial\chi_\alpha|^2e^{-\chi_\alpha}
 = (|f|^2+\alpha)^{s-3}
 (\alpha+s|f|^2)^2|\partial f|^2
 \le |f|^{p-2}|\partial f|^2,
\end{aligned}
\]
the dominated convergence theorem yields
\begin{equation}\label{eq:p-regularized-energy}
 \lim_{\alpha\downarrow0}\int_\Omega d\,
 |\partial f-f\partial\chi_\alpha|^2e^{-\chi_\alpha}\,dV
 =\frac14\int_\Omega d\,\Delta|f|^p\,dV
 \le C_\Omega\|f\|_{H^p}^p.
\end{equation}
For $1/2\le r<1$, put
\[
 E_{r,\alpha}:=\int_\Omega d^{2-r}
 \left(|v|^2+\langle[\Theta,\Lambda]^{-1}
 v_{g,\alpha},v_{g,\alpha}\rangle\right)e^{-\varphi_\alpha}\,dV.
\]
By \eqref{eq:p-regularized-energy},
$d^{2-r}\le C_\Omega d$, and \eqref{eq:unified-source-bound},
\begin{equation}\label{eq:Hardy-source-bound}
 \limsup_{\alpha\downarrow0}E_{r,\alpha}
 \le C_\Omega\delta^{-2q-2}|\log\delta|^2\|f\|_{H^p}^p,
\end{equation}
where $C_\Omega$ is independent of $p$, $m$, and $r$.

Fix $r$ and $\alpha>0$
with $E_{r,\alpha}<\infty$, and denote by
$E_{\varepsilon,r,\alpha}$ the same integral over
$\Omega_\varepsilon$, and set
$\operatorname{div}L=\sum_{\lambda=1}^n\partial_\lambda L_\lambda$.
Since $\rho_\varepsilon=0$ on
$\partial\Omega_\varepsilon$, $L\rho_\varepsilon=\vartheta$,
and $1-\vartheta$ vanishes near $\partial\Omega$,
\begin{equation}\label{eq:normal-pairing-coeff}
 |\rho_\varepsilon|
 +|1-L\rho_\varepsilon-\rho_\varepsilon\operatorname{div}L|
 \le C_\Omega d.
\end{equation}
For $\omega\in\mathcal D_1(\Omega_\varepsilon)$,
the holomorphy of $f,g$ gives
$v\in W^{1,2}_{(0,1)}(\Omega_\varepsilon)^m$; integration by parts
yields
\[
\begin{aligned}
 (v,\omega)_{\varphi_\alpha}
 =\int_{\Omega_\varepsilon}
 \Big\langle(1-L\rho_\varepsilon
 -\rho_\varepsilon\operatorname{div}L)v
 -\rho_\varepsilon v_{g,\alpha},\omega\Big\rangle
 e^{-\varphi_\alpha}\,dV
-\int_{\Omega_\varepsilon}\rho_\varepsilon
 \langle v,\bar L\omega\rangle e^{-\varphi_\alpha}\,dV.
\end{aligned}
\]
Since
$[\Theta,\Lambda]\ge I$, $|\bar L\omega|\le
C_\Omega|\bar\nabla\omega|$, and $-\rho\asymp d$, 
Cauchy--Schwarz, \eqref{eq:normal-pairing-coeff} and \eqref{eq:profile-applied}
yield
\begin{align}\label{eq:normal-pairing-Hp}
 |(v,\omega)_{\varphi_\alpha}|^2
 &\le C_\Omega E_{\varepsilon,r,\alpha}
 \int_{\Omega_\varepsilon}(-\rho)^r
 \bigl(\langle[\Theta,\Lambda]\omega,\omega\rangle
 +|\bar\nabla\omega|^2\bigr)e^{-\varphi_\alpha}\,dV\notag\\
 &\le\frac{C_\Omega E_{\varepsilon,r,\alpha}}{1-r}
 \left(\|T^*_{\mathcal H,\varphi_\alpha}\omega\|_{(-\rho)^r,\varphi_\alpha}^2
 +\|S_{\mathcal H}\omega\|_{(-\rho)^r,\varphi_\alpha}^2\right).
\end{align}
The estimate extends to $\mathcal G_{\varphi_\alpha}$ by
graph-norm density. Let $\Pi$ be the orthogonal projection
of $\mathcal H_1$ onto $\ker S_{\mathcal H}$. Since
$\operatorname{Ran}T_{\mathcal H}\subset\ker S_{\mathcal H}$
and $v\in\ker S_{\mathcal H}$, for every
$\omega\in\operatorname{Dom}T^*_{\mathcal H,\varphi_\alpha}$,
\[
 \Pi\omega\in\mathcal G_{\varphi_\alpha}\cap\ker S_{\mathcal H},
 \quad T^*_{\mathcal H,\varphi_\alpha}\Pi\omega
 =T^*_{\mathcal H,\varphi_\alpha}\omega,
 \quad (v,\Pi\omega)_{\varphi_\alpha}
 =(v,\omega)_{\varphi_\alpha}.
\]
Applying \eqref{eq:normal-pairing-Hp} to $\Pi\omega$,
\[
 |(v,\omega)_{\varphi_\alpha}|^2
 \le\frac{C_\Omega E_{\varepsilon,r,\alpha}}{1-r}
 \|T^*_{\mathcal H,\varphi_\alpha}\omega\|_{(-\rho)^r,\varphi_\alpha}^2.
\]
Hahn--Banach and Riesz representation give, after multiplication by $(-\rho)^r$,
$u_{\varepsilon,r,\alpha}$ satisfying
\[
\begin{gathered}
 (u_{\varepsilon,r,\alpha},
 T^*_{\mathcal H,\varphi_\alpha}\omega)_{\varphi_\alpha}
 =(v,\omega)_{\varphi_\alpha},\ \ \ 
 g\cdot u_{\varepsilon,r,\alpha}=0,\\
 (1-r)\int_{\Omega_\varepsilon}|u_{\varepsilon,r,\alpha}|^2
 (-\rho)^{-r}e^{-\varphi_\alpha}\,dV
 \le C_\Omega E_{\varepsilon,r,\alpha}.
\end{gathered}
\]
Hence $\dbar u_{\varepsilon,r,\alpha}=v$.
As $E_{\varepsilon,r,\alpha}\le E_{r,\alpha}$,
local weak $L^2$ compactness and lower semicontinuity as
$\varepsilon\downarrow0$ yield $u_{r,\alpha}$ on $\Omega$
with $\dbar u_{r,\alpha}=v$, $g\cdot u_{r,\alpha}=0$, and
\[
 (1-r)\int_\Omega|u_{r,\alpha}|^2
 (-\rho)^{-r}e^{-\varphi_\alpha}\,dV
 \le C_\Omega E_{r,\alpha}.
\]
Since $e^{-\varphi_\alpha}\ge c_\Omega e^{-\chi_\alpha}$
and $(-\rho)^{-r}\asymp d^{-r}$ uniformly for $1/2\le r<1$,
\eqref{eq:Hardy-source-bound} implies
\begin{equation}\label{eq:Hardy-correction-estimate}
 \limsup_{\alpha\downarrow0}(1-r)
 \int_\Omega|u_{r,\alpha}|^2e^{-\chi_\alpha}d^{-r}\,dV
 \le C_\Omega\delta^{-2q-2}
 |\log\delta|^2\|f\|_{H^p}^p.
\end{equation}
Now put
\[
 I_r=(1-r)\int_\Omega|f|^p d^{-r}\,dV.
\]
For each $c>0$, $e^{-\chi_c}\le e^{-\chi_\alpha}$ when
$0<\alpha<c$. Local weak $L^2$ compactness, lower
semicontinuity, and then monotone convergence as $c\downarrow0$
give, along a sequence $\alpha\downarrow0$, a limit $u_r$ with
\(
 \dbar u_r=v,\ g\cdot u_r=0,
\)
and
\[
 (1-r)\int_\Omega |u_r|^2|f|^{p-2}d^{-r}\,dV
 \le C_\Omega\delta^{-2q-2}|\log\delta|^2\|f\|_{H^p}^p.
\]

Set $h_r=f\bar g/|g|^2-u_r$. Then $\dbar h_r=0$ in
distributions.
Moreover, $g\cdot h_r=f$ and
\[
 |h_r|^2=\frac{|f|^2}{|g|^2}+|u_r|^2
\]
almost everywhere, since $g\cdot u_r=0$. Hence
\[
 (1-r)\int_\Omega |h_r|^2|f|^{p-2}d^{-r}\,dV
 \le\delta^{-2}I_r+
 C_\Omega\delta^{-2q-2}|\log\delta|^2\|f\|_{H^p}^p.
\]
For $p<2$, H\"older's inequality with exponents
$2/p$ and $2/(2-p)$ gives
\[
\begin{aligned}
\left((1-r)\int_\Omega|h_r|^p d^{-r}\,dV\right)^{1/p}
&\le
\left((1-r)\int_\Omega
|h_r|^2|f|^{p-2}d^{-r}\,dV\right)^{1/2}
I_r^{1/p-1/2}\\
&\le
\left(\delta^{-2}I_r+
C_\Omega\delta^{-2q-2}|\log\delta|^2
\|f\|_{H^p}^p\right)^{1/2}
I_r^{1/p-1/2}.
\end{aligned}
\]
For $p=2$, the same inequality follows directly.
Applying Lemma~\ref{lm:Hardy-limit-Hp} with a defining
function equal to $-d$ near $\partial\Omega$, we have
$I_r\to\|f\|_{H^p}^p$ as $r\uparrow1$, and a subsequence
of $h_r$ converges locally uniformly to $h\in H^p(\Omega)^m$.
Thus $g\cdot h=f$ and
\[
\begin{aligned}
 \|h\|_{H^p}
 \le\left(\delta^{-2}+
 C_\Omega\delta^{-2q-2}|\log\delta|^2\right)^{1/2}
 \|f\|_{H^p}
  \le C_\Omega\delta^{-q-1}|\log\delta|\,\|f\|_{H^p}.
\end{aligned}
\]
The estimates above are uniform in $p$, so the final constant
$C_\Omega$ is independent of $p$.
\end{proof}

\section{An example where the $H^p$-Corona Theorem fails for $p>2$}
\label{sec:counterexample}
Let $\varphi$ be a smooth nondecreasing convex function on $[0,\infty)$,
zero on $[0,1/16]$, with $\varphi(1/2)>4$ and constant positive
derivative on $[1,\infty)$. Put
\[
 \psi(t)=\int_0^t e^{-(16/y)^7}\,dy,\ \ \  t\ge0.
\]
$ \psi$ is convex and flat at zero, with $0<\psi'(t)\le1$ for $t>0$.
Let $k_0$ be a sufficiently large integer, to be chosen later.
For $k\ge k_0$ set $\xi_k=e^{i/k}$,
$a_k=e^{-k^{14}}$, and
\[
 Q_k(z_0)=e^{k^{12}(\overline{\xi_k}z_0-1)},\ \ \ 
 P_k(u,v)=(2u)^k\bigl(1+(2v)^{k^2}\bigr)-1.
\]
Define $\Omega=\{\rho<0\}\subset \C\times\C^2$ by
\begin{equation}\label{eq:ex-domain}
 \begin{aligned}
 \rho_0(z_0,z)&=|z_0|^2-1+\psi\bigl(|1-z_0|^2+|z|^2\bigr),\\
 \rho(z_0,z)&=\rho_0(z_0,z)
       +\sum_{k\ge k_0}a_k\varphi\bigl(|Q_k(z_0)P_k(kz)|^2\bigr).
 \end{aligned}
\end{equation}
Here $(z_0,z)\in\C\times\C^2$, with $z=(z_1,z_2)\in\C^2$. Thus

\begin{proposition}\label{ex:upper-range}
For $k_0$ sufficiently large, $\Omega$ is a bounded pseudoconvex domain with
$C^\infty$ boundary, which is strongly pseudoconvex
except at the infinite-type point $(1,0,0)$.
There is $g\in H^\infty(\Omega)^3$ satisfying
\eqref{eq:corona-data} for which $g\cdot h=1$ has no
solution in $H^p(\Omega)^3$ for any $p>2$.
\end{proposition}

\begin{proof}
Since $\varphi'$ is constant on $[1,\infty)$ and, for each compact set
$K\Subset\C^3$ and $j\ge0$,
\[
\|Q_k(z_0)P_k(kz)\|_{C^j(K)}
\le C_{K,j}e^{C_{K,j}k^{12}},
\]
Cauchy's estimates and the chain rule yield
\[
\sum_{k\ge k_0}
\left\|a_k\varphi\bigl(|Q_k(z_0)P_k(kz)|^2\bigr)
\right\|_{C^j(K)}
\le C_{K,j}\sum_{k\ge k_0}
e^{-k^{14}+C_{K,j}k^{12}}
\longrightarrow0.
\]
Thus $\rho$ is smooth and plurisubharmonic, and
${K_0}:=\{\rho_0\le0\}$ is compact with $\{\rho\le0\}\subset {K_0}$.
Since $\rho_0$ is convex and $\rho_0(0)=\psi(1)-1<0$, for
$(z_0,z)\in\{\rho=0\}$ we have
\[
\begin{aligned}
\left.\frac{d}{dr}\rho(rz_0,rz)\right|_{r=1}
&=\nabla\rho_0(z_0,z)\cdot(z_0,z)
 +\nabla(\rho-\rho_0)(z_0,z)\cdot(z_0,z)\\
&\ge 1-\psi(1)-(\rho-\rho_0)(z_0,z)
 -\sup_{K_0}|(z_0,z)|\,
 \|\rho-\rho_0\|_{C^1(K_0)}\\
&\ge 1-\psi(1)
 -\bigl(1+\sup_{K_0}|(z_0,z)|\bigr)
 \|\rho-\rho_0\|_{C^1(K_0)}
>0
\end{aligned}
\]
for sufficiently large $k_0$. Since $\rho(0)<0$, the radial derivative estimate also shows
that $\Omega$ is star-shaped. Thus $\Omega$ is a bounded
pseudoconvex domain with $C^\infty$ boundary and
plurisubharmonic defining function $\rho$.

Since $\psi'(t)>0$ for $t>0$, $\rho$ is strictly plurisubharmonic
away from $(1,0,0)$. At this point
\[
 P_k(0)=-1,\ \ \
 |Q_k(1)|=e^{-k^{12}(1-\cos(1/k))}<1/4
\]
for every $k\ge k_0$. Hence \( \forall\, k\ge k_0,\ \exists\,\varepsilon_k>0\) such that \(
\varphi\bigl(|Q_k(z_0)P_k(kz)|^2\bigr)
\equiv0
\ \text{on }B_{\varepsilon_k}(1,0,0),\)
\[
\begin{gathered}
D^\alpha(\rho-\rho_0)(1,0,0)=0,
\ \ \ |\alpha|\ge0.
\end{gathered}
\]
Since $\rho_0(1,\tau,0)=\psi(|\tau|^2)$ is flat at $\tau=0$,
$(1,0,0)\in\partial\Omega$ is the unique weakly pseudoconvex
boundary point and has infinite type.

Let $B$ be the Blaschke product in $\mathbb D$ with zeros
\[
x_{k,\ell}=(1-a_k)\xi_k e^{10i\ell a_k},
\ \ \ 
\ell\in\mathbb Z,\quad |10\ell a_k|\le k^{-8},
\]
namely, for $z_0\in\mathbb D,$
\[
B(z_0)=
\prod_{k\ge k_0}\prod_\ell
\frac{|x_{k,\ell}|}{x_{k,\ell}}
\frac{x_{k,\ell}-z_0}
{1-\overline{x_{k,\ell}}z_0}.
\]
Set
\[
J_k:=\#\{\ell\in\mathbb Z:|10\ell a_k|\le k^{-8}\}
=2\left\lfloor\frac{k^{-8}}{10a_k}\right\rfloor+1.
\]
Since
\[
\sum_{k\ge k_0}\sum_\ell(1-|x_{k,\ell}|)
=\sum_{k\ge k_0}J_ka_k
\asymp\sum_{k\ge k_0}k^{-8}<\infty,
\]
the corresponding Blaschke product $B$ is well defined. Set
\[
 Q(z_0)=\sum_{k\ge k_0} k Q_k(z_0)z_0^{L_k},\ \ \ 
 L_k=\lceil a_k^{-1}\rceil .
\]
The series converges locally uniformly in $\mathbb D$.
Set
\[
 g=(g_1,g_2,g_3):=c_0\bigl(Q(z_0)z_1,Q(z_0)z_2,B(z_0)\bigr),
\]
where $c_0>0$ is chosen so that $\sup_\Omega|g|=1$.
Since $Q$ and $1/B$
are bounded on $|z_0|\le1/2$, it suffices to consider $1/2<|z_0|<1$.
Put $t=1-|z_0|$ and let
\[
I_k:=\left\{\xi_k e^{i\theta}:
|\theta|<\frac{1}{10k^2}\right\}
\subset\partial\mathbb D.
\]
These arcs are disjoint, and
$|Q_k(z_0)|\le e^{-ck^8}$ whenever $z_0/|z_0|\notin I_k$. Thus
\[
\begin{gathered}
\sum_{\substack{j\ge k_0\\ z_0/|z_0|\notin I_j}}
j|Q_j(z_0)|\,|z_0|^{L_j}
\le \sum_{j\ge k_0}je^{-cj^8}<\infty,\\
k|Q_k(z_0)|\,|z_0|^{L_k}
\le ke^{-t/a_k},
\ \ \ \frac{z_0}{|z_0|}\in I_k.
\end{gathered}
\]

If $t>\sqrt{a_k}$, then
\(
ke^{-t/a_k}
\le ke^{-1/\sqrt{a_k}}
=ke^{-e^{k^{14}/2}}\le 6
\)
for sufficiently large $k$.
If $t\le\sqrt{a_k}$, $\psi(|z|^2)\le \psi(|1-z_0|^2+|z|^2)\le 1-|z_0|^2\le 2\sqrt{a_k}$ on $\Omega$.
Since $\psi(y)\ge(y/2)e^{-(32/y)^7}$, if $|z|\ge6/k$,
then
\[
2e^{-k^{14}/2}
\ge\psi(|z|^2)
\ge\frac{18}{k^2}
e^{-(8/9)^7k^{14}},
\]
which is impossible for sufficiently large $k$, since
$(8/9)^7<1/2$. Hence $k|z|<6$. Consequently,
 $|Q(z_0)|\,|z|$ is bounded
on $\Omega$.

For $x=x_{k,\ell}$ write $b_x(z_0)=(z_0-x)/(1-\bar x z_0)$.
For $z_0\ne0$, set
$\theta_{k,\ell}
=|\operatorname{Arg}(z_0\overline{x_{k,\ell}})|$.
Then
\[
1-|b_{x_{k,\ell}}(z_0)|^2
\le
\frac{Ca_k(1-|z_0|)}
{(a_k+1-|z_0|)^2+\theta_{k,\ell}^2}.
\]
If $z_0/|z_0|\notin I_k$, then
$\theta_{k,\ell}\gtrsim k^{-2}$, so
\[
\begin{aligned}
\sum_{\substack{k\ge k_0\\
z_0/|z_0|\notin I_k}}
\sum_\ell
\bigl(1-|b_{x_{k,\ell}}(z_0)|^2\bigr)
\le C\sum_{k\ge k_0}J_ka_k k^4
\lesssim\sum_{k\ge k_0}k^{-4}<\infty.
\end{aligned}
\]
If $z_0/|z_0|\in I_k$, since the arcs $I_k$ are pairwise disjoint,
\[
\#\left\{k\ge k_0:
\frac{z_0}{|z_0|}\in I_k\right\}\le1.
\]
For the possible remaining index $k$, set
$t=1-|z_0|$ and $A=a_k+t$.
Since \( \arg x_{k,\ell+1}-\arg x_{k,\ell}=10a_k\), choose
$\ell_0$ so that
$\theta_{k,\ell_0}=\min_\ell\theta_{k,\ell}$.
Then
\[
\theta_{k,\ell_0\pm j}
\ge (j-\tfrac12)10a_k,\ \ \  j\ge1,
\]
whenever the corresponding indices exist.
 The monotonicity
of $(A^2+y^2)^{-1}$ for $y\ge0$ yields
\[
\begin{aligned}
\sum_\ell\frac{1}{A^2+\theta_{k,\ell}^2}
\le \frac{1}{A^2}
 +2\sum_{j=1}^\infty
 \frac{1}{A^2+(j-\frac12)^2(10a_k)^2}
\le \frac{1}{A^2}
 +\frac{C}{a_k}
 \int_0^\infty\frac{dy}{A^2+y^2}.
\end{aligned}
\]
Consequently,
\[
\begin{aligned}
\sum_\ell
\bigl(1-|b_{x_{k,\ell}}(z_0)|^2\bigr)
\le C\left(
\frac{a_kt}{(a_k+t)^2}
+\int_{\mathbb R}
\frac{tdy}{(a_k+t)^2+y^2}
\right)
=C\left(
\frac{a_kt}{(a_k+t)^2}
+\frac{\pi t}{a_k+t}
\right)
\le C.
\end{aligned}
\]
Combining the two cases, we obtain
\[
\begin{aligned}
\sum_{k\ge k_0}\sum_\ell
\bigl(1-|b_{x_{k,\ell}}(z_0)|^2\bigr)
\le C\sum_{k\ge k_0}k^{-4}
+C\,\#\left\{k\ge k_0:
\frac{z_0}{|z_0|}\in I_k\right\}
\le C.
\end{aligned}
\]
Since $-\log y\le C(1-y^2)$ on $[1/10,1]$, $|B|$ is bounded below
if $|b_x(z_0)|\ge 1/10$. If $|b_x(z_0)|<1/10$, we have $|z_0-x|\le2a_k/9$,
$1-|z_0|^2<3a_k$, $|z_0|^{L_k}\ge c$, and
\[
 k^{12}\bigl(1-\operatorname{Re}(\overline{\xi_k}z_0)\bigr)
 =O\bigl(k^{12}(a_k+k^{-16})\bigr)=o(1).
\]
Consequently $|Q(z_0)|\asymp k$. If also $|kz|\le1/4$, then
$P_k(kz)=-1+O(2^{-k})$, so
$|Q_k(z_0)P_k(kz)|^2>1/2$ for large $k$.
Thus $a_k\varphi(|Q_k(z_0)P_k(kz)|^2)>4a_k$,
while $\rho_0\ge-3a_k$, so $\rho(z_0,z)>a_k>0$.
Hence
\[
\{(z_0,z)\in\C^3:
|b_{x_{k,\ell}}(z_0)|<1/10,\ |kz|\le1/4\}
\cap\Omega=\varnothing.
\]
Therefore, when $|b_{x_{k,\ell}}(z_0)|<1/10$
and $(z_0,z)\in\Omega$, we have
\[
|Q(z_0)|\,|z|\ge ck|z|>c/4.
\]
Consequently,
$g$ satisfies \eqref{eq:corona-data}.

Fix $k\ge k_0$ and an admissible index $\ell$. We claim that
\begin{equation}\label{eq:ex-collar}
 |z_0-x_{k,\ell}|<a_k/8,\
 |(kz_1,kz_2)|<2,\
 |P_k(kz_1,kz_2)|<1/4
 \ \Longrightarrow\
 (z_0,z)\in\Omega\cap\{d<2a_k\}.
\end{equation}
Indeed, $|1-z_0|<2/k$ and $|z|<2/k$, so
\[
 \psi\bigl(|1-z_0|^2+|z|^2\bigr)
 \le \psi(8/k^2)
 \le (8/k^2)e^{-128k^{14}}=o(a_k).
\]
Since $|x_{k,\ell}|=1-a_k$ and
$|z_0-x_{k,\ell}|<a_k/8$,
\[
1-|z_0|^2\ge 1-|z_0|\ge\frac78a_k,
\ \ \ 
|Q_k(z_0)P_k(kz_1,kz_2)|<\frac14.
\]
Moreover, as $\arg\xi_k=1/k$,
\[
\begin{aligned}
\left|\arg z_0-\frac1k\right|
\le|\arg z_0-\arg x_{k,\ell}|
+\left|\arg x_{k,\ell}-\frac1k\right|
\le Ca_k+k^{-8}<\frac1{10k^2}
\end{aligned}
\]
for sufficiently large $k$. Hence
$z_0/|z_0|\in I_k$, and therefore
$z_0/|z_0|\notin I_j$ for every $j\ne k$.
Thus, for $j\ne k$,
\[
1-\operatorname{Re}(\overline{\xi_j}z_0)
\ge cj^{-4},
\ \ \ 
|P_j(jz_1,jz_2)|\le e^{Cj^2\log j},
\]
and
\[
|Q_j(z_0)P_j(jz_1,jz_2)|
\le e^{-cj^8+Cj^2\log j}<\frac14
\ \ \ (j\ge k_0,\ j\ne k).
\]
Since $\varphi=0$ on $[0,1/16]$, every
summand in $\rho-\rho_0$ vanishes. Consequently,
\[
\rho(z_0,z)=\rho_0(z_0,z)
\le-\frac78a_k+\psi(8/k^2)<0.
\]
Since $\Omega\subset\{|z_0|<1\}$,
\[
d(z_0,z)\le1-|z_0|
\le\frac98a_k<2a_k.
\]

Suppose there exists $h\in H^p(\Omega)^3$ such that
$g\cdot h=1$. On the slice $\{z_0=x_{k,\ell}\}$,
we have $B(x_{k,\ell})=0$, so
\[
kz_1h_1(x_{k,\ell},z_1,z_2)
+kz_2h_2(x_{k,\ell},z_1,z_2)
=\frac{k}{c_0Q(x_{k,\ell})}.
\]
Following Oh's modification of Cole's construction
\cite[pp.~300--301]{Oh}, for each \(\zeta\) with \(2^{-k}\le|\zeta|\le2^{k^2}\), consider the \(k^3\) solutions, counted with multiplicity, of
\[
(2kz_1)^k=\frac1\zeta,\ \ \ 
(2kz_2)^{k^2}=\zeta-1.
\]
They satisfy $P_k(kz_1,kz_2)=0$ and
$|(kz_1,kz_2)|<2$.
Let $A(\zeta)$ be the average of
$kz_1h_1(x_{k,\ell},z_1,z_2)$ over these roots,
counted with multiplicity.
Since analytic continuation permutes the roots,
$A$ is single-valued and holomorphic on
$\{2^{-k}<|\zeta|<2^{k^2}\}\setminus\{1\}$.
As $\zeta\to1$, the roots remain in a compact subset
of the slice of $\Omega$, so $A$ is locally bounded
and extends across $\zeta=1$ by the removable
singularity theorem. Since $z_2=0$ at $\zeta=1$,
\[
A(1)=\frac{k}{c_0Q(x_{k,\ell})},
\ \ \  |A(1)|\asymp1.
\]
On $|\zeta|=2^{k^2}$, we have $|kz_1|=2^{-k-1}$,
and on $|\zeta|=2^{-k}$, we have $|kz_1|=1$.
Let $M_{k,\ell}$ be the maximum of $|h_1|$
over the roots above these two circles. For $\zeta$
on either circle, the triangle inequality gives
\[
\begin{aligned}
|A(\zeta)|
\le\frac1{k^3}\sum_{j=1}^{k^3}
|kz_1^{(j)}|\,
|h_1(x_{k,\ell},z_1^{(j)},z_2^{(j)})|
\le\frac12|\zeta|^{-1/k}M_{k,\ell}.
\end{aligned}
\]
Then
\[
\max_{|\zeta|=2^{k^2}}|A(\zeta)|
\le 2^{-k-1}M_{k,\ell},
\ \ \ 
\max_{|\zeta|=2^{-k}}|A(\zeta)|
\le M_{k,\ell}.
\]
Cauchy's formula on the annulus
$\{2^{-k}<|\zeta|<2^{k^2}\}$ yields
\[
\begin{aligned}
c\le |A(1)|
\le\left(
2^{-k-1}\frac{2^{k^2}}{2^{k^2}-1}
+\frac{2^{-k}}{1-2^{-k}}
\right)M_{k,\ell}
\lesssim2^{-k}M_{k,\ell}.
\end{aligned}
\]

For each $(k,\ell)$, choose $\zeta_*$ on one of
the two circles and a corresponding root
$(z_{1,*},z_{2,*})$ such that
\[
|h_1(x_{k,\ell},z_{1,*},z_{2,*})|
=M_{k,\ell}.
\]
Consider the bidisk
\[
\Delta_{k,\ell}:=\bigl\{
(z_{1,*}(1+\alpha),z_{2,*}(1+\beta)):
|\alpha|<k^{-3},\quad
|\beta|<k^{-3}\min\{1,|\zeta_*|\}
\bigr\}.
\]
It lies in $|(kz_1,kz_2)|<2$. Since
$(2kz_{1,*})^k=\zeta_*^{-1}$ and
$(2kz_{2,*})^{k^2}=\zeta_*-1$, throughout
$\Delta_{k,\ell}$
\[
\begin{aligned}
P_k(kz_1,kz_2)
=(1+\alpha)^k-1
+(1+\alpha)^k(1-\zeta_*^{-1})
\bigl((1+\beta)^{k^2}-1\bigr)
=O(k^{-1}).
\end{aligned}
\]
Indeed, $(1+\alpha)^k-1=O(k^{-2})$ and
$(1+\beta)^{k^2}-1=O(k^2|\beta|)$.
Since
\(
|kz_{1,*}|=1/2|\zeta_*|^{-1/k}
\) and \(
|kz_{2,*}|=1/2|\zeta_*-1|^{1/k^2},
\)
we have
\[
|kz_{1,*}|\min\{1,|\zeta_*|\}
=
\begin{cases}
2^{-k}, & |\zeta_*|=2^{-k},\\
2^{-k-1}, & |\zeta_*|=2^{k^2}.
\end{cases}
\]
Moreover,
\[
|\zeta_*-1|\asymp
\begin{cases}
1, & |\zeta_*|=2^{-k},\\
2^{k^2}, & |\zeta_*|=2^{k^2}.
\end{cases}
\]
Consequently,
\(
|kz_{1,*}kz_{2,*}|^2
\min\{1,|\zeta_*|\}^2\asymp2^{-2k}.
\)
Thus the polydisk
\[
\mathcal P_{k,\ell}
:=\{z_0\in\mathbb C:|z_0-x_{k,\ell}|<a_k/8\}
\times\Delta_{k,\ell}
\] satisfies
\[
\begin{aligned}
\operatorname{Vol}(\mathcal P_{k,\ell})
=\frac{\pi^3a_k^2}{64k^{16}}
|kz_{1,*}\,kz_{2,*}|^2
\min\{1,|\zeta_*|\}^2
\asymp a_k^2k^{-16}2^{-2k}.
\end{aligned}
\]
For fixed $k$, these polydisks $\mathcal P_{k,\ell}$ are pairwise disjoint
as $\ell$ varies. Moreover, by \eqref{eq:ex-collar},
\[
\mathcal P_{k,\ell}
\subset\Omega\cap\{d<2a_k\}.
\]
The submean inequality and the disjointness of
these polydisks give
\[
\begin{aligned}
\int_{\Omega\cap\{d<2a_k\}}|h_1|^p\,dV
&\ge\sum_\ell
\int_{\mathcal P_{k,\ell}}|h_1|^p\,dV\\
&\ge\sum_\ell
\operatorname{Vol}(\mathcal P_{k,\ell})
\left|h_1(x_{k,\ell},
z_{1,*},z_{2,*})\right|^p\\
&=\sum_\ell
\operatorname{Vol}(\mathcal P_{k,\ell})M_{k,\ell}^p\\
&\ge c_pJ_ka_k^2k^{-16}2^{(p-2)k}.
\end{aligned}
\]
Applying the coarea formula to the boundary
distance $d$, and \eqref{eq:Hardy-means-Hp}
with a defining function equal to $-d$ near
$\partial\Omega$, we obtain
\[
\begin{aligned}
\int_{\Omega\cap\{d<2a_k\}}|h_1|^p\,dV
=\int_0^{2a_k}
\int_{\{d=s\}}|h_1|^p\,d\sigma_s\,ds
\le C_\Omega a_k\|h_1\|_{H^p(\Omega)}^p.
\end{aligned}
\]
Since $J_ka_k\asymp k^{-8}$,
\[
\begin{aligned}
c_p'a_kk^{-24}2^{(p-2)k}
\le c_pJ_ka_k^2k^{-16}2^{(p-2)k}
\le\int_{\Omega\cap\{d<2a_k\}}|h_1|^p\,dV
\le C_\Omega a_k\|h_1\|_{H^p(\Omega)}^p.
\end{aligned}
\]
Consequently,
\[
k^{-24}2^{(p-2)k}
\le C_p\|h_1\|_{H^p(\Omega)}^p,
\]
which is impossible for $p>2$ as $k\to\infty$.
\end{proof}

\section*{Declaration of AI assistance}
The author used ChatGPT (OpenAI) to assist with the development and checking
of proofs and with writing, revision, and proofreading. In particular,
ChatGPT devised the construction of the counterexample in
Proposition~\ref{ex:upper-range}. The author reviewed
and edited all AI-generated content and takes full responsibility
for the content of this article.

\end{document}